\documentclass[10pt]{article}

\usepackage{amsmath,amssymb,latexsym}
\usepackage{array,graphicx,xcolor}
\usepackage{hyperref}

\newtheorem{theo}{Theorem}[section]

\newtheorem{prop}[theo]{Proposition}
\newtheorem{rem}[theo]{Remark}

\newtheorem{lemm}[theo]{Lemma}

\def\endproof{\hfill{}$\square$}
\def\la{\lambda}

\newcommand{\proof}{\medskip\goodbreak\noindent {\textbf{Proof.} }}
\newcommand{\dis}{\displaystyle}

\newcommand\RR{\mathbb{R}}

\newcommand\sig{\sigma}
\newcommand\ds{\displaystyle}

\providecommand{\abs}[1]{\lvert#1\rvert}
\providecommand{\norm}[1]{\lVert#1\rVert}

\title{Lipschitz stability in an inverse problem for the Kuramoto-Sivashinsky equation}
\author{Lucie Baudouin\footnotemark[1], 
 \,Eduardo Cerpa\footnotemark[2], 
 \,Emmanuelle Cr\'epeau\footnotemark[3], 
 \,and \,Alberto Mercado\footnotemark[2]}
 
\begin{document}
 \date{}

\maketitle
\footnotetext[1]{CNRS, LAAS, 7 avenue du colonel Roche, F-31400 Toulouse, France.\\
Univ de Toulouse, LAAS, F-31400 Toulouse, France.\\
 E-mail: {\tt lucie.baudouin@laas.fr}}
\footnotetext[2]{Departamento de Matem\'atica, Universidad T\'ecnica Federico Santa Mar\'ia, Casilla 110-V, Valpara\'iso, Chile.\\
 E-mail: {\tt eduardo.cerpa@usm.cl, alberto.mercado@usm.cl}}
\footnotetext[3]{Laboratoire de Math\'ematiques, Universit\'e de Versailles Saint-Quentin en Yvelines, 78035 Versailles, France.\\
 E-mail: {\tt emmanuelle.crepeau@math.uvsq.fr}}

\abstract{ In this article, we present an inverse problem for the nonlinear 1-d Kuramoto-Sivashinsky (K-S) equation.
More precisely, we study the nonlinear inverse problem of retrieving the anti-diffusion coefficient from the measurements of the solution on a part of the boundary and also at some positive time in the whole space domain. The Lipschitz stability for this inverse problem is our main result and it relies on the Bukhge{\u\i}m-Klibanov method. The proof is indeed based on a global Carleman estimate for the linearized K-S equation.}

\vspace{0.2 cm}

{\bf Keywords:} Inverse problem, Kuramoto-Sivashinsky equation, Carleman estimate

\vspace{0.2 cm} {\bf AMS subject classifications:} 35R30, 35K55.

\section{Introduction}\label{intro}

We focus in this paper on an inverse problem that consists in the determination of a coe\-ffi\-cient in a partial differential equation (pde) from the partial knowledge of a given single solution of the equation. This class of problems (single-measurement 
{ coefficient} 
inverse problems) was investigated using Carleman estimates  for the first time  in \cite{81bukhgeim-klibanov}  
by Bukhge{\u\i}m and Klibanov. 
{ 
See \cite{Klibanov92}, \cite{91klibanov-malinsky} 
and the recent book \cite{Beilina-Klibanov} for details about the  so-called  Bukhge{\u\i}m-Klibanov method}. 
This method {  was initially used} to prove uniqueness for inverse problems (i.e. that each measurement corresponds to only one coefficient) from local Carleman estimates
{ (estimates valid for solutions with compact support in the interior of the domain), as in \cite{81bukhgeim-klibanov}.}
Regarding the continuity { of the inverse problem of recovering the source term,} the first Lipschitz stability result for a multidimensional wave equation was obtained by Puel and Yamamoto \cite{96puel-yamamoto} 
{ using the uniqueness result and a compactness-uniqueness argument.
}

{  Global Carleman estimates (valid for solutions considered in the whole domain and satisfying boundary conditions) were applied to parabolic equations by first time in \cite{98imanuvilov-yamamoto}, where Lipschitz stability of an inverse problem is established.} 
  Since then, this type of inverse problems for parabolic equations has received a large amount of attention. The primary difference with respect to hyperbolic inverse problems is that parabolic problems are not time-reversible:  therefore, an additional measurement must be added if that method is applied. As one can read in the discussion of the introduction of \cite{98imanuvilov-yamamoto}, the knowledge of the full-state of the solution for some positive time is required. 
Proving the Lipschitz stability without  this assumption, which  is usually needed when global Carleman inequalities are used, is still an open problem.
Nevertheless, there are some uniqueness results with less assumptions on the measurements, that can be found in the litterature, such as \cite{10cristofol-roques} or some other inversion method in \cite{Klibanov92,KliTiBook}.

Recent results regarding linear parabolic problems can be found in \cite{07benabdallah-gaitan-lerousseau} (discontinuous coefficient),  \cite{06cristofol-gaitan-ramoul} (systems), \cite{Ignat-Pazoto-Rosier} (network) and the references therein. In \cite{09boulakia-grandmont-osses, 05egger-engl-klibanov,  10cristofol-roques}, nonlinear parabolic equations were even considered.  

{  
Among others pde's  coefficient/source inverse problems where Carleman estimates have been used   we can mention, without being exhaustive, logarithmic stability \cite{06bellassoued-yamamoto}, Calder\'on problem \cite{KenigUhlmann2009} or Schr\"odinger equation \cite{02baudouin-puel}. }  

In this paper, we consider a 1D non\-linear fourth-order parabolic equation called Kuramoto-Sivashinsky (K-S) equation. 
This equation was proposed independently by  Kuramoto and Tsuzuki \cite{Kur-Tzu} as a model for the phase turbulence in reaction diffusion systems, and by Sivashinsky \cite{Siv}, as a model the physical phenomena of plane flame propagation, were
the combined influence of diffusion and thermal conduction of a gas is described.\\

The K-S equation with non-constant coefficients describing the diffusion $\sigma=\sigma(x)$, and the anti-diffusion $\gamma=\gamma(x)$, is given as \begin{equation}
\left\{\begin{array}{ll}\label{KS1}
y_t+(\sigma(x)y_{xx})_{xx}+\gamma(x) y_{xx}+yy_x=g, &\qquad \forall (t,x)\in Q,\\
y(t,0)=h_1(t),\quad y(t,1)=h_2(t),&\qquad \forall t\in(0,T),\\
y_x(t,0)=h_3(t), \quad y_x(t,1)=h_4(t),&\qquad \forall t\in(0,T),\\
y(0,x)= y_0(x),&\qquad \forall x\in(0,1),
\end{array}\right.
\end{equation}
where $Q:=(0,T)\times(0,1)$,  $\sigma:[0,1]\rightarrow\RR^*_+$, and the functions $y_0, g, h_j$ are the initial condition,  the source term 
and the boundary data respectively. All these terms are assumed to be known and compatible.

In this nonlinear pde, the fourth-order term models the diffusion, and the second-order term models the incipient instabilities. 
We consider the inverse problem of retrieving the anti-diffusion coefficient $\gamma$ from boun\-da\-ry measurements of the solution.  This corresponds for instance  to getting information on the instability of a reaction-diffusion media  by measuring a single solution, which could represent a flame propagating on the domain. { Concerning the boundary measurements we will make, it is worth to mention that in a fourth-order parabolic problem like KS, boundary data $u_{xx}$ and $u_{xxx}$ are referred to as Neumann data, which in fact represent heat flux \cite{liu-krstic} in this kind of models}.\\

To the knowledge of the authors there are no results in the literature concerning the determination of coefficients for this nonlinear equation. However, a Carleman estimate has been used to obtain the null-controllability of the K-S equation in reference \cite{cerpa-mercado} for the constant coefficient case. { Other results on the control of the KS equation can be found in \cite{hu-temam,armou,armou2,liu-krstic,10cerpa}}.\\

Since the linearized equation is parabolic, we know that boundary measurements will not be sufficient to prove stability and we must consider an additional measurement of the full solution for a given time $T_0$ (as in \cite{07benabdallah-gaitan-lerousseau, 98imanuvilov-yamamoto} among others).\\

Our first result involves the local well-posedness of the nonlinear equation \eqref{KS1}. A less regular framework can be used for this equation { but} the method applied in this paper requires the solution and its time-derivative to be at least in $L^2(0,T;H^4(0,1))$. Therefore, let us introduce the following notations for the functional spaces appearing in this paper:
\begin{equation}\label{notation} 
\begin{array}{l} 
\mathcal Y_k:= C([0,T]; H^k(0,1))\cap L^2(0,T;H^{k+2}(0,1)), \quad \text{ for } k\in\mathbb N; \\
\mathcal F:= \{f\in L^2(0,T;H^4(0,1)) \big / \, f_t\in L^2(0,T;L^2(0,1))\};\\
\mathcal Z := \{z\in \mathcal Y_6\big / \,z_t\in \mathcal Y_2\}.
\end{array} 
\end{equation}

\begin{theo}\label{wp} 
Let $\gamma\in H^4(0,1)$ and  $\sigma \in H^4(0,1)$ be such that 
\begin{equation}\label{sig}
\sigma(x)\geq \sigma_0>0,~~ \forall x\in (0,1).
\end{equation}
There exists $\varepsilon>0$ such that if  $y_0 \in H^6(0,1)$,  $g\in \mathcal F$, and $h_j \in H^2(0,T)$ for $j=1, \ldots, 4$
satisfy the compatibility conditions 
\begin{equation}\label{compcond}
y_0(0) = h_1(0),\quad y_{0,x}(0) = h_3(0), \quad y_0(1) = h_2(0), \quad y_{0,x}(1) = h_4(0),
\end{equation}
and
\begin{equation}\label{epsilon}
\|y_0\|_{H^6(0,1)}\leq \varepsilon,\quad \|g\|_{\mathcal F}\leq \varepsilon, \quad \|h_j\|_{H^2(0,T)}\leq \varepsilon \, \mbox{ for  }  \,  j=1, \ldots, 4,
\end{equation} 
then the K-S equation \eqref{KS1}  has a unique solution $y\in \mathcal Z$ .
\end{theo}

Once the existence of solutions to the K-S equation has been established (see Section~\ref{well}), the following inverse problem is addressed:
\begin{quote}
Is it possible to retrieve the anti-diffusion coefficient $\gamma=\gamma(x)$ from the measurement of $y_{xx}(t,0)$ and $y_{xxx}(t,0)$ on $(0,T)$ and from the measurement of $y(T_0,x)$ on $(0,1)$, where $y$ is the solution to Equation (\ref{KS1}) and $T_0 \in (0,T)$?
\end{quote}
A local answer for this nonlinear inverse problem is given (see section~\ref{InversePb}). To be more specific, let  $\tilde\gamma$ be given and  fixed. We denote by $\tilde y$ the solution to Equation (\ref{KS1}) with $\gamma$ replaced by $\tilde\gamma$. This paper focuses on the following question concerning the unknown $\gamma$ and $y$.

\noindent\textbf{Stability}: Is it possible to estimate $\|\tilde \gamma - \gamma\|_{L^2(0,1)}$ by suitable norms $\|\tilde y (T_0,x) - y (T_0,x)\|$ in space and $\|\tilde y_{xx}(t,0) - y_{xx}(t,0)\|$, $\|\tilde y_{xxx} (t,0)- y_{xxx}(t,0)\|$ in time?

Of course, a positive answer implies the usual uniqueness result.

\noindent\textbf{Uniqueness}: Do the equalities of the measurements $\tilde y_{xx}(t,0) = y_{xx}(t,0)$ and $\tilde y_{xxx}(t,0) = y_{xxx}(t,0)$ for $t\in(0,T)$ and $\tilde y (T_0,x) = y (T_0,x) $ for $x\in(0,1)$ imply $\tilde\gamma = \gamma$ on $(0,1)$?\\

In order to answer these questions, we use the Bukhge{\u\i}m-Klibanov method. First, a global Carleman estimate for the linearized K-S equation with non-constant coefficients is obtained. It is then used to prove the main result which can be stated as follows.

To precisely state the results we prove in this article, we introduce, for $m > 0$, the set
$$
    L^\infty_{\leq m} (0,1) =  \left\{ \gamma \in L^\infty(0,1) s.t. \  \norm{\gamma}_{L^\infty(0,1)} \leq m  \right\}.
$$

\begin{theo}\label{PIKS} Let us consider $\sigma\in H^4(0,1)$ satisfying \eqref{sig}, $\gamma \in H^{4}(0,1)$, $g\in \mathcal F$ and the data $y_0 \in H^6(0,1)$ and $h_j \in H^2(0,T)$ for $j=1, \ldots, 4$ under the compatibility conditions \eqref{compcond}. 
Let  $y \in \mathcal Z$ be the solution of \eqref{KS1},  and $\tilde y \in \mathcal Z$ the solution corresponding to a given $\tilde\gamma \in H^4(0,1)$ instead of~$\gamma$. 
We assume that there exists $\eta>0$ and  $T_0\in(0,T)$ such that
\begin{equation}
	\label{InitialDataCond}
	 \inf \left\{\left|\tilde y_{xx}(T_0,x)\right|,x\in (0,1)\right\} \geq \eta ,
\end{equation}

Then, given $M > 0$, there exists a positive constant $C$ depending on the parameters $(T, m,M ,\eta)$, such that for every $\gamma\in  L^\infty_{\leq m} (0,1)$,
\begin{multline}\label{stabNL}
	\|\gamma-\tilde\gamma\|^2_{L^2(0,1)} \leq C\left\| y_{xx}(\cdot,0) - \tilde y_{xx}(\cdot,0)\right\|^2_{H^1(0,T)}
	+ C\left\| y_{xxx}(\cdot,0) - \tilde y_{xxx}(\cdot,0)\right\|^2_{H^1(0,T)}\\
	+ C\left\|y (T_0,\cdot) - \tilde y (T_0,\cdot) \right\|^2_{H^4(0,1)}+ C\left\|y (T_0,\cdot) - \tilde y (T_0,\cdot) \right\|^4_{H^1(0,1)}
\end{multline}
for all $y$ satisfying
$$
	\|y\|_{\mathcal Z}\leq M.
$$
\end{theo}

This inequality states the stability of the inverse problem. Before giving the outline of our paper and the proofs of the different steps, we want to give several comments on this result.

\begin{rem}
{ For numerical purposes it would be interesting to know explicitly how the constant $C$ in \eqref{stabNL} depends on the diffusion $\sigma$ or on the time $T$. This kind of question has been addressed in \cite{coron-guerrero,fernandez-cara-guerrero} for observability constant in the framework of second-order parabolic equations. In those papers the authors got an exponential dependence on both the constant diffusion and the time.}
\end{rem}

\begin{rem} 
One can show that there exist solutions satisfying assumption \eqref{InitialDataCond}. We present two different arguments:
\begin{enumerate}
\item We take $\varepsilon > 0$ given  by Theorem \ref{wp}, and some $y^0 \in H^6(0,1)$ such that $\ds{\inf_{x \in (0,1)}  \left| y^0_{xx} \right | \geq \varepsilon/2}$. For arbitrary boundary data and source term belonging to the corresponding spaces, by Theorem \ref{wp} there exists a solution  $\tilde y \in C([0,1];H^6(0,1))$ with $\tilde y(0,\cdot) = y^0$. Using Sobolev injection and continuity, we obtain the existence of a time $T_0 >0$ such that \eqref{InitialDataCond} is fulfilled with $\eta = \varepsilon/4$.
\item 
We can also prove that there exist solutions satisfying  \eqref{InitialDataCond} without asking $T_0$ to be small, but instead, constraining the source term and boundary data as follows:  Let $y_0$ be the initial data and $T_0$ belong to $(0,T)$. Let us  pick up a state $y_1=y_1(x)$ strictly convex. We consider the trajectory $\tilde y(t,x)=\frac{T_0 -t}{T_0}y_0(x) + \frac{t}{T_0}y_1(x)$, which is the solution of equation \eqref{KS1} with  source term given by $g=\tilde y_t+(\sigma(x)\tilde y_{xx})_{xx}+\gamma(x) \tilde y_{xx}+\tilde y\tilde y_x$ and the boundary data  given by the traces of $\tilde y$. Thus,  $\tilde y(T_0,x)=y_1(x)$ and hence the trajectory $\tilde y$ satisfies \eqref{InitialDataCond}.
\end{enumerate}
Therefore, the set of data and solutions where our stability result is valid is not empty.
\end{rem}

\begin{rem} 
We obtain the same result if $\tilde y$ has a different initial condition than $y$. See in Section 4, that the term $v(x,0)$ of system \eqref{KS3} does no play any role in the result.
\end{rem}

\begin{rem}
We can complete inequality \eqref{stabNL} by the following:
\begin{multline*}
	\left\| y_{xx}(\cdot,0) - \tilde y_{xx}(\cdot,0)\right\|^2_{H^1(0,T)}
	+ \left\| y_{xxx}(\cdot,0) - \tilde y_{xxx}(\cdot,0)\right\|^2_{H^1(0,T)}\\
	+ \left\|y (T_0,\cdot) - \tilde y (T_0,\cdot) \right\|^2_{H^4(0,1)}+ \left\|y (T_0,\cdot) - \tilde y (T_0,\cdot) \right\|^4_{H^1(0,1)}\\
	\leq C\Big(\|y-\tilde y\|^2_{H^1(0,T;H^4(0,1))}+\|y-\tilde y\|^4_{C([0,T];H^1(0,1))}\Big).
\end{multline*}
This inequality follows directly from standard Sobolev injections. It  indicates that the required measurements are finite if $y$ and $\tilde y$ belong to the space $H^1(0,T;H^4(0,1))$ and this is true if $y$ and $\tilde y$ are solutions in $\mathcal Z$ provided by Theorem \ref{wp}. 
\end{rem}

\begin{rem} As stated in the introduction, an internal measurement at $t=T_0$ is required if this method, using Carleman estimates, is used to prove the stability for this type of inverse problem for parabolic equations. Nevertheless, this is probably a technical point since there is no counter-example that demonstrates whether this assumption is required for stability. In $\cite{10cristofol-roques}$, uniqueness (but not stability) is proven using a very different technique in an inverse problem for a parabolic equation and without any internal measurements in the whole space domain. One can also mention a method in \cite{KliTiBook} that can deliver uniqueness from hyperbolic equations to parabolic ones.
\end{rem}

\begin{rem} 
In this paper, the boundary measurements are located at $x=0$, but the result would be the same if we measure at $x=1$ instead. Indeed, the choice of a suitable weight function in the proof of the Carleman estimate in Section~$\ref{CarlemanEstimate}$ is critical to impose the side of measurement.
\end{rem}

This article is organized as follows. The well-posedness result stated in Theorem~\ref{wp} is proved in  Section~\ref{well}.  A global Carleman estimate for a general K-S equation is given and proved in Section~\ref{CarlemanEstimate}. Finally, Section~\ref{InversePb} contains the use of the Bukhgeim-Klibanov method to prove the Lipschitz stability of the inverse problem stated in Theorem \ref{PIKS}.

\section{ On the Cauchy problem for KS equation}\label{well}

This section presents a proof of Theorem \ref{wp}  in a more general case including time dependent lower-order coefficients. 
We consider the following K-S system
\begin{equation}
	\left\{\begin{array}{ll}\label{KSt}
		y_t+(\sigma(x)y_{xx})_{xx}+\gamma(x) y_{xx}+G_1 y_x + G_2 y+ yy_x=g, & \qquad\forall (t,x)\in Q,\\
		y(t,0)=h_1(t), \quad y(t,1)=h_2(t),& \qquad\forall t\in(0,T),\\
		y_x(t,0)=h_3(t), \quad y_x(t,1)=h_4(t),& \qquad\forall t\in(0,T),\\
		y(0,x)= y_0(x),& \qquad\forall x\in(0,1),
	\end{array}\right.
\end{equation}
where $G_1,G_2$ belong to $H^1(0,T; H^4(0,1))$, $g\in\mathcal F$ and $y_0\in H^6(0,1)$ is compatible with $h_j \in H^2(0,T)$ for $j=1, \ldots , 4$. 
Recall that the coefficients satisfy $\gamma\in H^4(0,1)$, $\sigma\in H^4(0,1)$ and hypothesis \eqref{sig}.  \\

First, we only consider the main part of the linear differential operator in the next proposition.

\begin{prop}\label{p1} Let $z_0\in H^6\cap H_0^2(0,1)$ and $f\in \mathcal F$. Then, the following equation
\begin{equation}
	\left\{\begin{array}{ll}\label{KS-P}
		 z_t+(\sigma(x) z_{xx})_{xx}=f, & \qquad\forall (t,x)\in Q,\\
		z(t,0)=0,\quad z(t,1)=0,&\qquad \forall t\in (0,T),\\
		z_x(t,0)=0, \quad z_x(t,1)=0,&\qquad \forall t\in (0,T),\\
		z(0,x)= z_0(x), &\qquad \forall x\in (0,1),
	\end{array}\right.
\end{equation}
has a unique solution $z\in \mathcal Z$ and there exists $C>0$ such that   
$$ 
	\| z \|_\mathcal Z \leq C \left(  \| f \|_{\mathcal F} + \|z_0\|_{H^6}  \right).
$$
\end{prop}

\proof The operator 
$${  	
	\begin{array}{rcl}
		H^4\cap H_0^2(0,1)\subset L^2(0,1)  &\longrightarrow&  L^2(0,1)\\
	z & \longmapsto& (\sigma(x) z''(x))'' ,
	\end{array}
	}
$$
is simultaneously positive, coercive and self-adjoint. 
By the Hille-Yosida-Phillips Theorem (see \cite{Caz-Har}),  it generates a strongly continuous semigroup in $L^2(0,1)$. Therefore, for each  $z_0\in H^4\cap H^2_0(0,1)$ and $f\in C^1([0,T];L^2(0,1))$, Equation  \eqref{KS-P} 
has a unique solution $z \in C([0,T];H^4\cap H_0^2(0,1)) \cap C^1([0,T];L^2(0,1))$. 

We will demonstrate that the solutions $z \in \mathcal Z$ (refer to the notation introduced in \eqref{notation}), can be obtained by taking $z_0$ and $f$ sufficiently regular.
 
We now search for some energy estimates that indicate the space where the solutions lie on depending on the regularity of the data. Suppose that there are solutions sufficiently regular to perform the following computations. Equation \eqref{KS-P} is multiplied by $z$ and integrated over $(0,1)$ in space. 
Some integrations by parts give
\begin{equation}
\label{e2}
\frac{d}{dt}\left(\int_0^1 |z(t,x)|^2\,dx \right)+\int_0^1|z_{xx}(t,x)|^2\,dx\leq C \left( \int_0^1 |f(t,x)|^2\,dx+ \int_0^1 |z(t,x)|^2\,dx \right).
\end{equation}
Throughout this paper, $C$ denotes a positive constant that may vary from line to line. To make the reading easier, we denote for any function $u$ of $x$ and $t$, $$\iint_Q u= \int_0^T\int_0^1 u(t,x) \,dxdt$$

Using Gronwall's lemma, we first obtain that for all $t>0$,
\begin{equation}
\label{energy1}
\int_0^1 |z(t,x)|^2\,dx\leq C \left( \iint_Q |f|^2+ \int_0^1 |z_0|^2\,dx \right).
\end{equation}

\noindent Then, (\ref{e2}) is integrated over $[0,T]$ and (\ref{energy1}) is used to get
 \begin{equation}
\label{energy2}
\iint_Q |z_{xx}|^2 \leq C \left(\iint_Q |f|^2+  \int_0^1 |z_0|^2\,dx\right).
\end{equation}
Inequalities (\ref{energy1}) and (\ref{energy2}) finally imply that
\begin{equation}\label{Bprima}\|z\|^2_{\mathcal Y_0}\leq C \iint_Q |f|^2+C \int_0^1 |z_0|^2\,dx .\end{equation}

Now, equation \eqref{KS-P} is multiplied by $(\sigma z_{xx})_{xx}$ and integrated over $(0,1)$ in space. 
Some integrations by parts give also
$$ 
\frac{1}{2}\frac{d}{dt}\left(\int_0^1\sigma | z_{xx}(t,x)|^2\,dx\right) + \int_0^1|(\sigma z_{xx}(t,x))_{xx}|^2\,dx=\int_0^1f(t,x)(\sigma z_{xx}(t,x))_{xx}\,dx.
$$

\noindent Using the inequality  $ab\leq \frac 12 a^2 + \frac 12 b^2$, we get
\begin{equation} \label{Eqsigma1}
\frac{d}{dt}\left(\int_0^1 \sigma| z_{xx}(t,x)|^2\,dx\right)+ \int_0^1|( \sigma  z_{xx}(t,x))_{xx}|^2\,dx \leq \int_0^1 |f(t,x)|^2\,dx.
\end{equation}
Using Gronwall's lemma, from\eqref{Eqsigma1} and \eqref{sig} we obtain that for all $t>0$,
\begin{equation}
\label{energy10}
\int_0^1 |z_{xx}(t,x)|^2\,dx\leq C \left( \iint_Q |f|^2+ \int_0^1 |z_0''|^2\,dx \right).
\end{equation}

\noindent Then, \eqref{Eqsigma1} is integrated over $[0,T]$ and (\ref{energy10}) is used to get
 \begin{equation}
\label{energy20}
\iint_Q |(\sigma z_{xx})_{xx}|^2 \leq C \left(\iint_Q |f|^2+  \int_0^1 |z_0''|^2\,dx\right),
\end{equation}
and then, taking into account that $\sigma \in H^4$, we get
 \begin{equation}
\label{energy21}
\iint_Q | z_{xxxx}|^2 \leq C \left(\iint_Q |f|^2+  \int_0^1 |z_0''|^2\,dx\right) + C \|z\|_{L^2(0,T;H^3(0,1))}.
\end{equation}
For any $\varepsilon > 0$, from Ehrling's Lemma (see Theorem 7.30 in \cite{ReRo}) and \eqref{energy1}, we have that 
 \begin{equation} \label{energy22}
 	\begin{aligned}
		\ds{ \|z\|_{L^2(0,T;H^3(0,1))} } & \ds{ \leq 
		\varepsilon \|z\|_{L^2(0,T;H^4(0,1))} + C \|z\|_{L^2(0,T;L^2(0,1))} }\\
		& \leq  \varepsilon \|z\|_{L^2(0,T;H^4(0,1))} + C \left( \iint_Q |f|^2+ \int_0^1 |z_0|^2\,dx \right).
	 \end{aligned}
\end{equation}
Taking $\varepsilon > 0$ small enough, inequalities \eqref{energy10},  \eqref{energy21}  and  \eqref{energy22}  imply that
\begin{equation}\label{e}
	\|z\|^2_{\mathcal Y_2}\leq C\iint_Q |f|^2+ C \|z_0\|_{L^2(0,T;H^2(0,1))} ^2.
\end{equation}

On the other hand, Equation \eqref{KS-P} is derived with respect to time. Thus $q:=z_t$ satisfies
\begin{equation}
	\left\{\begin{array}{ll}\label{KS-T}
		 q_t+(\sigma(x) q_{xx})_{xx}=f_t, & \qquad\forall (t,x)\in Q,\\
		q(t,0)=0,\quad q(t,1)=0,&\qquad \forall t\in (0,T),\\
		q_x(t,0)=0, \quad q_x(t,1)=0,&\qquad \forall t\in (0,T),\\
		q(0,x)= f(0,x)-(\sigma z^{\prime\prime}_0(x))^{\prime\prime}, &\qquad \forall x\in (0,1).
	\end{array}\right.
\end{equation}

Using estimate \eqref{e}, we obtain $q\in \mathcal Y_2$ if $(f(0,x)-(\sigma z^{\prime\prime}_0(x))^{\prime\prime})\in H^2(0,1)$ and $f_t\in L^2(0,T;L^2(0,1))$. These hypotheses are fulfilled if $z_0\in H^6\cap H_0^2(0,1)$ and $f\in \mathcal F $. Note that $\mathcal F\subset C([0,T];H^2(0,1)).$ From the equation satisfied by $z$ and the fact that $f\in\mathcal F$ and $z_t\in \mathcal Y_2$, we determine that $z\in \mathcal  Y_6$, which concludes the proof of Proposition \ref{p1}. \endproof \\

Then, we focus on  the linear problem with non-homogenous boundary conditions and  low-order coefficients that depend on time.
\begin{prop} 
Let $z_0\in H^6(0,1)$, $\hat f \in \mathcal F$,  $G_1,G_2\in H^1(0,T; H^4(0,1))$
and $h_j \in H^2(0,T)$ for $j=1, \ldots , 4$ satisfying the compatibility conditions with $z_0$. 
Then, the equation
\begin{equation}
	\left\{\begin{array}{ll}\label{KS-NL}
		 z_t+(\sigma(x) z_{xx})_{xx} +\gamma(x) z_{xx}+G_1z_x+G_2 z =\hat f, & \qquad\forall (t,x)\in Q,\\
		z(t,0)=h_1(t),\quad z(t,1)=h_2(t),&\qquad \forall t\in (0,T),\\
		z_x(t,0)=h_3(t), \quad z_x(t,1)=h_4(t),&\qquad \forall t\in (0,T),\\
		z(0,x)= z_0(x), &\qquad\forall x\in (0,1),
	\end{array}\right.
\end{equation}
has a unique solution $z\in \mathcal Z$ and there exists $C>0$ such that   
$$ 
	\| z \|_\mathcal Z \leq C \left(  \|Ê\hat f \|_{\mathcal F} + \|z_0\|_{H^6}  + \sum_{j=1}^{4} \|h_j\|_{H^2} \right).
$$
\end{prop}
\proof 
We first prove this result for null boundary data (i.e. for $h_j = 0$ for $j=1, \ldots , 4$ and therefore $z_0\in H^6\cap H^2_0(0,1)$). \\
For any $\hat w\in \mathcal Z$, $\Pi(\hat w)$ is defined as the solution of  \eqref{KS-P} with $f=(\hat f-\gamma(x) \hat w_{xx}-G_1\hat w_x-G_2 \hat w)$. Note that $f\in \mathcal F$ and therefore $\Pi(\hat w)\in \mathcal Z$ is well defined.

If $T$ is small enough, then $\Pi$ is a contraction. Indeed, for any $w,\hat w\in \mathcal Z$, we have
\begin{eqnarray}\|\Pi(\hat w)-\Pi(w)\|_{\mathcal Z} &\leq&C \|\gamma(x)(w_{xx}- \hat w_{xx})+G_1(w_x-\hat w_x)+G_2 (w-\hat w)\|_{\mathcal F} \nonumber\\
&\leq&C\|w-\hat w\|_{L^2(H^6)}+C\|w_t-\hat w_t\|_{L^2(H^2)}\label{m1}\\
&\leq&CT^{\frac 1 4}\|w-\hat w\|_{L^4(H^6)}+CT^{\frac 1 4}\|w_t-\hat w_t\|_{L^4(H^2)} \nonumber\\
&\leq&CT^{\frac 1 4}\|w-\hat w\|_{\mathcal Y_6}+CT^{\frac 1 4}\|w_t-\hat w_t\|_{\mathcal Y_2} \nonumber\\
&\leq&CT^{\frac 1 4}\|w-\hat w\|_{\mathcal Z}\label{m4},
 \end{eqnarray}
where the space $L^m(0,T;H^n(0,1))$ is denoted as $L^m(H^n)$.
 
Hence, the operator $\Pi$ has a unique fixed point in $\mathcal Z$, which is the solution of \eqref{KS-NL} with $h_j = 0$ for $j=1, \ldots , 4$. Using standard arguments and the linearity of this equation, the solution can be extended to a larger time interval.

In order to prove the general case, take $h_j \in H^2(0,T)$,  $j=1, \ldots , 4$ compatible with~$z_0$.
It is not difficult to find a function $\psi \in H^2(0,T; C^\infty([0,1]))$ satisfying the boundary conditions of \eqref{KS-NL}. 
For instance take $\psi(x,t) = \sum_{j=1}^4p_j(x)h_j(t)$ 
where $p_1(x) = 2x^3 - 3x^2 +1$, $p_2(x) = -2x^3 + 3x^2$, $p_3(x) = x^3 - 2x^2 +x$ and $p_4(x) = x^3 - x^2$.
In particular we have $L\psi :=   \psi_t+(\sigma(x)  \psi_{xx})_{xx} +\gamma(x)  \psi_{xx}+G_1 \psi_x+G_2 \psi  \in \cal{F}$. Then, if $w$ is the solution of equation \eqref{KS-NL}  with null boundary data,
initial condition $w_0  - \psi(\cdot,0)$,
and right-hand side equal to $\hat f  - L\psi$,
 let us define $z = w + \psi$.
It is not difficult to see that $z$ is the required solution.
\quad \endproof

\begin{rem} The third-order term $z_{xxx}$ can be added to Equation \eqref{KS-NL}. Indeed, in that case \eqref{m1} becomes $C\|w-\hat w\|_{L^2(H^7)}+C\|w_t-\hat w_t\|_{L^2(H^3)}$, which is bounded by $$CT^{\frac 1 4}\|w-\hat w\|^{1/2}_{L^\infty(H^6)}\|w-\hat w\|^{1/2}_{L^2(H^8)}+CT^{\frac 1 4}\|w_t-\hat w_t\|^{1/2}_{L^\infty(H^2)}\|w_t-\hat w_t\|^{1/2}_{L^2(H^4)}.$$ This last expression is bounded by \eqref{m4}. The remainder of the proof is the same.
\end{rem}

Again, by using a fixed point theorem, we can prove Theorem \ref{wp} for equation \eqref{KSt}.\\
Let $y_0 \in H^6(0,1)$, $h_j\in H^2(0,1)$ compatible with $y_0$, and $g\in \mathcal F$.
For any $v\in \mathcal Z$, we define $\Lambda(v)$ as the solution of  \eqref{KS-NL} with $\hat f=( g-vv_x)$ and $z_0=y_0$. Note that $\hat f \in \mathcal F$ and therefore $\Lambda(v)\in \mathcal Z$ is well defined. Indeed, if $v\in\mathcal Y_3$ and $v_t\in \mathcal Y_0$, then we have
$$(v v_x)_{xxxx}=(10v_{xx}v_{xxx}+5v_{x}v_{xxxx}+v v_{xxxxx})\in L^2(0,T;L^2(0,1))$$ and
$$(vv_x)_t=v_tv_x+vv_{xt}\in L^2(0,T;L^2(0,1)).$$
Furthermore, we can prove 
\begin{equation}\label{rr}
\begin{aligned}
\|\Lambda(v)\|_{\mathcal Z}  & \leq C \left(\|g\|_{\mathcal F}+\|vv_x\|_{\mathcal F}+\|y_0\|_{H^6}  + \sum_{j=1}^4\|h_j\|_{H^2} \right) \\
					    & \leq C \left(\|g\|_{\mathcal F}+\|v\|^2_{\mathcal Z}+\|y_0\|_{H^6} + \sum_{j=1}^4\|h_j\|_{H^2} \right).
\end{aligned}
\end{equation}

Let $\varepsilon>0$ 
and suppose that $y_0$, $h_j$ and $g$ satisfy \eqref{epsilon}. Consider $v$ such that $\|v\|_{\mathcal Z}\leq R$ with $R>0$ satisfying 
$C(6\varepsilon + R^2)<R$. 
From \eqref{rr}, we obtain $\|\Lambda(v)\|_{\mathcal Z}<R$. 
Thus, the application $\Lambda$  maps the ball $B_R:=\{v\in \mathcal Z \big/ \|v\|_{\mathcal Z}\leq R\} $ into itself.

 We will now prove that $\Lambda:B_R \rightarrow B_R$ is a contraction. 
 For any $z,v\in B_R$,  $\Lambda(z)-\Lambda(v)$ is the solution of \eqref{KS-NL} with 
 $z_0=0$, $h_j=0$ for $j=1, \ldots, 4$ and 
 $\hat f = vv_x - zz_x$. We obtain the estimate
 $$
 	\|\Lambda(z)-\Lambda(v)\|_{\mathcal Z}
 		\leq C \|vv_x-zz_x\|_{\mathcal F}
		\leq C\left(\|(v-z)v_x\|_{\mathcal F} + \| z(v_x-z_x)\|_{\mathcal F}\right).
$$
Using the definition \eqref{notation} of the space $\mathcal F$, 
$v,z\in C([0,1];H^6(0,1))   \hookrightarrow L^\infty(0,T;W^{5,\infty}(0,1))$ and 
$v_t,z_t\in C([0,1];H^6(0,1))  \hookrightarrow  L^\infty(0,T;W^{1,\infty}(0,1))$, we obtain
$$
	\|\Lambda(z)-\Lambda(v)\|_{\mathcal Z}
	\leq C(\|v\|_{\mathcal Z} + \|z\|_{\mathcal Z})\| v-z\|_{\mathcal Z}
	\leq 2CR\| v-z\|_{\mathcal Z},
$$ 
which implies that $\Lambda$ is a contraction if $R$ is chosen small enough. More precisely, we can choose $R,\varepsilon$ such that $2CR<1$ and $C(6\varepsilon+R^2)<R$. Hence, the map $\Lambda$ has a unique fixed point $y\in\mathcal Z$, which is the unique solution of \eqref{KSt}. 
Thus, we have proven Theorem \ref{wp}. \endproof

\section{Global Carleman inequality}\label{CarlemanEstimate}

In this section, a global Carleman inequality will be proved for the linearized K-S equation. 
We define the space
\begin{equation}\label{HIP}
\mathcal V = \{ v \in L^2(0,T; H^4\cap H_0^2(0,1)) \, \,  \big/ \, \,   Lv \in L^2((0,T)\times(0,1)) \}
\end{equation}
 where 
\begin{equation*}
Lv=v_t+(\sigma v_{xx})_{xx}+q_2v_{xx}+q_1v_x+q_0v
\end{equation*}
with  $q_j \in L^\infty(\Omega)$ for $j = 0,1,2$.\\
Consider $\beta \in C^4([0,1])$  such that for some $r > 0$ we have, for all $x\in (0,1)$:
\begin{equation}\label{hip1B}
\begin{array}{l}
0 <  r  \leq  \beta (x),  \\
0 <  r  \leq  \beta' (x),  \\
\beta''(x) \leq - r < 0,   \\
\abs{\sig'(x)\beta'(x)}\leq \dfrac{r}{4} \dis\min_{ z \in[0,1]}\{\sig(z)\}.  
\end{array}
\end{equation}
For instance, if $\sigma$ is constant, we can consider $\beta(x) =  \sqrt{1+x}$.\\
On the other hand, given $T_0 \in (0,T)$  we can choose $\phi_0 \in C^1([0,T])$ such that
\begin{equation}\label{hip1P}
\begin{array}{l}
\phi_0(0) = \phi_0(T) = 0,  \, \, \mbox{ and }\\
0 < \phi_0(t) \leq  \phi_0(T_0) \, \,  \mbox{ for each } t \in (0,T).  
\end{array}
\end{equation}
For example, if $T_0 = T/2$, we can use $\phi_0(t) = t(T-t)$.\\

We finally define the function 
\begin{equation}\label{defphi}
\phi(t,x) = \frac{\beta(x)}{\phi_0(t)}, 
\end{equation}
for $(t,x) \in (0,T) \times [0,1]$,
which is the weight function of the Carleman estimate. 
From \eqref{hip1B}  and \eqref{hip1P} it is not difficult to see that $\phi$ satisfies the following properties:
\begin{equation} \label{Phi1}
\begin{array}{l}
\exists C>0 \mbox{ such that } \phi \leq C \phi_x \mbox{ and } \\
\phi^n \leq C \phi^{m} \, \,  \mbox{ for each positive integers } n < m.
\end{array}
\end{equation}

\begin{theo}\label{CarlemanKS} Let $\phi$ be a function defined by $(\ref{defphi})$
and $m > 0$.
Then there exists $\lambda_0>0$ and a constant $C = C(T, \lambda_0, r,m)>0$ such that  
if
 $\|q_i\|_{L^\infty((0,T) \times (0,1))}\leq m$ for $i=0,1,2$
then  we have 
\begin{multline}\label{Carleman}
	\int_0^T\int_0^1 e^{-2\lambda \phi}\left(\frac{\abs{v_t}^2+\abs{(\sig v_{xx})_{xx}}^2}{\lambda\phi}
	+\lambda^7\phi^7\abs{v}^2+\lambda^5\phi^5\abs{v_x}^2+\lambda^3\phi^3\abs{v_{xx}}^2+\lambda\phi\abs{v_{xxx}}^2\right)dxdt\\
	\leq C \int_0^T\int_0^1 e^{-2\lambda\phi}\abs{Lv}^2\,dxdt\\
		+ C\int_0^T e^{-2\lambda\phi(t,0)} \bigg(\lambda^3\phi_x^3(t,0)\sig(0)^2\abs{v_{xx}(t,0)}^2
		+\lambda\phi_x(t,0)\sig^2(0)\abs{v_{xxx}(t,0)}^2\bigg)\,dt
\end{multline}
for all $v \in\mathcal V$, for all $\lambda\geq \lambda_0$. 
\end{theo}

As we pointed out in the Introduction, a Carleman estimate for the K-S equation with constant coefficients $\sigma$ and $\gamma$  was previously obteined in  \cite{cerpa-mercado}. The final goal in that work was to prove null-controllability with boundary controls. Thus, \eqref{Carleman} is a generalization  to the case of non-constant coefficients.

\proof
Consider the following operator $P$ defined in $\mathcal W_\lambda := \{ e^{-\lambda \phi}v \, : \,  v  \in\mathcal V\}$ by
\begin{equation*}
Pw=e^{-\lambda \phi}L(e^{\lambda \phi}w).
\end{equation*}
We then obtain the decomposition  $Pw=P_1w+P_2w+Rw$, where
\begin{eqnarray}
 P_1w &=&6\lambda^2\phi_x^2\sig w_{xx}+\lambda^4\phi_x^4\sig w+(\sig w_{xx})_{xx}+6\lambda^2(\phi_x^2\sig)_xw_x \label{decomposition P1}
\\
P_2w&=&w_t+4\lambda^3\phi_x^3\sig w_x+4\lambda \phi_x\sig w_{xxx}+4\lambda^3\phi_x(\phi_x^2\sig)_xw\label{decomposition P2}\\
Rw&=&\lambda\phi_tw+2\lambda\phi_x\sig_{xx}w_x+\lambda^2\phi_x^2\sig_{xx}w+\lambda\phi_{xx}\sig_{xx}w\nonumber\\
&& + \,6\lambda\phi_x\sig_xw_{xx}+6\lambda^2\phi_x\phi_{xx}\sig_xw+
6\lambda\phi_{xx}\sig_xw+2\lambda\phi_{xxx}\sig_xw \nonumber\\
&& +\,  4\lambda^2\phi_x\phi_{xxx}\sig w+6\lambda\phi_{xx}\sig w_{xx}+3\lambda^2\phi_{xx}^2\sig w+4\lambda \phi_{xxx}\sig w_x\nonumber\\
&& + \, \lambda\phi_{xxxx}\sig w+q_0w+q_1w_x+q_1\lambda\phi_x w\nonumber\\
&& + \,q_2w_{xx}+2\lambda q_2\phi_xw_x+\lambda^2q_2\phi_x^2w+\lambda\phi_{xx}q_2w \nonumber\\
&& -  \,2\lambda^3\phi_x^2\phi_{xx}\sig w-2\lambda^3\phi_x^3\sig_x w. \label{defR}
\end{eqnarray}
Thus, $$ \|Pw  - Rw\|^2_{L^2(Q)} = \|P_1w\|_{L^2(Q)}^2 +  2\left< P_1w, P_2w\right> + \|P_2w\|_{L^2(Q)}^2 $$
where $\left< \cdot, \cdot \right>$ is the $L^2(Q)$ scalar product.

For any  $v \in \mathcal V$ we obtain $v_t\in L^2(0,T; L^2(0,1))$ and then $v\in C([0,T];L^2(0,1))$. From the construction of $\phi$ (see \eqref{hip1P}), we obtain $w\in C([0,T];L^2(0,1))$ and  $w(x,0) = w(x,T)=0$ for any $w \in \mathcal W_\lambda$.\\

Let us define the notations 
$$I(w)=-6\lambda^7\int_0^T\int_0^1\phi_x^6\phi_{xx}\sig^2\abs{w}^2\,dxdt,$$
$$I(w_x)=-\lambda^5\int_0^T\int_0^1\phi_x^4\sig(30\phi_{xx}\sig+12\phi_x\sig_x)\abs{w_x}^2\,dxdt,$$
$$I(w_{2x})=-\lambda^3\int_0^T\int_0^1\phi_x^2\sig(58\phi_{xx}\sig+40\phi_x\sig_x)\abs{w_{xx}}^2\,dxdt,$$
$$I(w_{3x})=-\lambda\int_0^T\int_0^1\sig(2\phi_{xx}\sig-4\phi_x\sig_x)\abs{w_{xxx}}^2\,dxdt,$$ and
$$I_x = \int_0^T (10 \lambda^3\phi_x^3\sig^2\abs{w_{xx}}^2 
         +  2 \lambda\phi_x\sig\sig_{xx}\abs{w_{xx}}^2 
         +  2  \lambda\phi_x\sig^2\abs{w_{xxx}}^2 ) \bigg|_{x=0}^1\,dt
$$
The following weighted norm is defined, for any $w \in \mathcal W_\lambda$, as
\begin{equation*}
 \left\| w \right\|_{_{\la, \phi}}^2 =
\int_0^T\int_0^1 \left(
\lambda^7\phi^7\abs{w}^2+\lambda^5\phi^5\abs{w_x}^2+\lambda^3\phi^3\abs{w_{xx}}^2+\lambda\phi\abs{w_{xxx}}^2\right)  \,dxdt. 
\end{equation*}

We first require the following 
\begin{lemm} \label{Lemma_prodL2}
Under  the  hypothesis of  Theorem~\ref{CarlemanKS}, there exists $\delta > 0$ such that 
\begin{equation} \label{prodL2}
\left<  P_1w, P_2w \right>_{L^2(Q)} \geq  \delta  \left\| w \right\|_{_{\la, \phi}}^2
+ I_x
\end{equation}
for $\lambda$ large enough and for all $w \in\mathcal W_\lambda$. 
\end{lemm}
\proof 
It is sufficient to prove that
\begin{equation} \label{prodL2-2}
\left<  P_1w, P_2w \right>_{L^2} =   \sum_{k=0}^3 I(w_{kx}) 
+ R_0(w) 
+ I_x
\end{equation}
for a large enough $\lambda$, for all $w \in\mathcal W_\lambda$, where $|R_0(w)| \leq   \lambda^{-1} \left\| w \right\|_{_{\la, \phi}}^2 $.

Indeed, let us first assume that we have (\ref{prodL2-2}). From the hypotheses in \eqref{hip1B} we easily check that there exists $\varepsilon>0$ such that $\phi$ satisfies for all $x \in (0,1)$,
\begin{equation}\label{desigphi1}
\begin{array}{ll}
\phi_{xx}(x)\leq - \varepsilon \phi < 0, \\
30\phi_{xx}(x)\sig(x)+12\phi_x(x)\sig_x(x)\leq - \varepsilon  \phi <  0,\\
58\phi_{xx}(x)\sig(x)+40\phi_x(x)\sig_x(x)\leq - \varepsilon  \phi <  0, \, \mbox{ and }\\
2\phi_{xx}(x)\sig(x)-4\phi_x(x)\sig_x(x)\leq - \varepsilon  \phi < 0. 
\end{array}
\end{equation}
Then from \eqref{Phi1} and assuming (\ref{prodL2-2}) we obtain, for $\lambda$ large enough,
\begin{equation} \label{prodL2-3}
\begin{aligned}
\left<  P_1w, P_2w \right>_{L^2} = &  \sum_{k=0}^3 I(w_{kx}) 
+ R_0(w) + I_x \\
& \geq  2\delta  \left\| w \right\|_{_{\la, \phi}}^2 - | R_0(w)|+ I_x \\
& \geq \delta  \left\| w \right\|_{_{\la, \phi}}^2 + I_x.
\end{aligned}
\end{equation}

Let us now prove (\ref{prodL2-2}): we write 
$\ds \left<  P_1w, P_2w \right>_{L^2(Q)} = \sum_{i,j=1}^4I_{i,j}$ 
where $I_{i,j}$ denotes the $L^2$-product between  the $i$-th term of $P_1w$ in \eqref{decomposition P1} and the $j$-th term of $P_2w$ in \eqref{decomposition P2}.\\
Integrations by parts in time or space are performed on each expression $I_{i,j}$. Each resulting expression will be included in one of the terms of the right-hand side of (\ref{prodL2-2}). The results are listed below, and we indicate for each term where it will be included.

\begin{itemize}
\item

$\ds{
I_{1,1}= - I_{4,1}
+
\underbrace{3\lambda^2\iint_Q(\phi_x^2\sig)_t\abs{w_x}^2}_{R_0(w)}}$
\item
$\ds{I_{1,2}=\underbrace{-12\lambda^5\iint_Q(\phi_x^5\sig^2)_x\abs{w_x}^2}_{I(w_x)}}$.
\item
$\ds{
I_{1,3}=\underbrace{-12\lambda^3\iint_Q (\phi_x^3\sig^2)_x\abs{w_{xx}}^2}_{I(w_{2x})}+\underbrace{12\lambda^3\int_0^T \bigg[\phi_x^3\sig^2\abs{w_{xx}}^2\bigg]_0^1 \,dt}_{I_x}}$
\item
$\ds{
I_{1,4}=\underbrace{12\lambda^5\iint_Q [\phi_x^3\sig(\phi_x^2\sig)_x]_{xx}\abs{w}^2}_{R_0(w)}
-\underbrace{24\lambda^5\iint_Q\phi_x^3\sig(\phi_x^2\sig)_x\abs{w_x}^2}_{I(w_x)}}$.
\item
$\ds{ I_{2,1}= \underbrace{ -\frac{\lambda^4}{2}\iint_Q(\phi_x^4\sig)_t\abs{w}^2}_{R_0(w)} }$.
\item
$\ds{
I_{2,2}=\underbrace{-2\lambda^7\iint_Q(\phi_x^7\sig^2)_x\abs{w}^2}_{I(w)}.
}$
\item
$\ds{
I_{2,3}=\underbrace{-2\lambda^5\iint_Q(\phi_x^5\sig^2)_{xxx}\abs{w}^2}_{R_0(w)}
+\underbrace{6\lambda^5\iint_Q(\phi_x^5\sig^2)_x\abs{w_x}^2}_{I(w_x)}}$.
\item
$\ds{
I_{2,4}=\underbrace{4\lambda^7\iint_Q \phi_x^5\sig(\phi_x^2\sig)_x\abs{w}^2}_{I(w)}
}$.
\item
$\ds{
I_{3,1}=  \frac{1}{2}\int_0^1\bigg[ \sig\abs{w_{xx}}^2\bigg]_0^T \,dx= 0}$.
\item
$\ds{
I_{3,2} = 
\underbrace{-2\lambda^3\iint_Q [(\phi_x^3\sig)_{xx}\sig]_x\abs{w_x}^2}_{R_0(w)}
+\underbrace{4\lambda^3\iint_Q(\phi_x^3\sig)_x\sig\abs{w_{xx}}^2}_{I(w_{2x})} } \\
~\hfill{}\ds{
+\underbrace{2\lambda^3\iint_Q(\phi_x^3)_x\sig^2\abs{w_{xx}}^2}_{I(w_{2x})}
-\underbrace{2\lambda^3\int_0^T\bigg[\phi_x^3\sig^2\abs{w_{xx}}^2\bigg]_0^1\,dt}_{I_x}
}$.
\item
$\ds{
I_{3,3}=\underbrace{2\lambda\int_0^T\bigg[\phi_x\sig\sig_{xx}\abs{w_{xx}}^2\bigg]_0^1\,dt}_{I_x}
-\underbrace{2\lambda\iint_Q (\phi_x\sig\sig_{xx})_x\abs{w_{xx}}^2}_{R_0(w)}
+\underbrace{8\lambda\iint_Q\phi_x\sig\sig_x\abs{w_{3x}}^2}_{I(w_{xxx})} }\\
~\hfill{}\ds{
+\underbrace{2\lambda\int_0^T\bigg[ \phi_x\sig^2\abs{w_{xxx}}^2\bigg]_0^1\,dt}_{I_x}
-\underbrace{2\lambda\iint_Q (\phi_x\sig^2)_x\abs{w_{xxx}}^2}_{I(w_{xxx})}}$.
\item
$\ds{
I_{3,4} =
\underbrace{4\lambda^3\iint_Q(\phi_x(\phi_x^2\sig)_x)_{xx}\sig ww_{xx}}_{R_0(w)}
-\underbrace{4\lambda^3\iint_Q (\phi_x(\phi_x^2\sig)_x)_{x}\sig\abs{w_x}^2}_{R_0(w)}}\\
~\hfill{}\ds{
+\underbrace{4\lambda^3\iint_Q \phi_x(\phi_x^2\sig)_x\sig\abs{w_{2x}}^2}_{I(w_{2x})}
}$.
\item 
$\ds{ 
I_{4,1}=6\lambda^2\iint_Q(\phi_x^2\sig)_xw_xw_t
}$, which is canceled when adding  with $I_{1,1}$.
\item
$\ds{
I_{4,2}=\underbrace{24\lambda^5\iint_Q (\phi_x^2\sig)_x\phi_x^3\sig\abs{w_x}^2}_{I(w_x)}
}$.
\item
$\ds{
I_{4,3}=\underbrace{12\lambda^3\iint_Q[(\phi_x^2\sig)_x\phi_x\sig]_{xx}\abs{w_x}^2}_{R_0(w)}
-\underbrace{24\lambda^3\iint_Q(\phi_x^2\sig)_x\phi_x\sig\abs{w_{xx}}^2}_{I(w_{2x})}
}$.
\item
$\ds{
I_{4,4}= \underbrace{-12\lambda^5\iint_Q (\phi_x^2\sig)_x(\phi_x^3\sig)_x\abs{w}^2}_{R_0(w)}
}$.
\end{itemize}
Summing up all the terms, we obtain (\ref{prodL2-2}).
\endproof\\

Then, we will  prove a Carleman inequality for the conjugated operator $P$.

\begin{lemm} \label{Carlemanw}
There exists $\lambda_0 > 0$ such that for all $\lambda \geq \lambda_0$ we have, for all $w \in\mathcal W_\lambda$,
\begin{eqnarray*}
 \int_0^T\int_0^1 \left(
\lambda^7\phi^7\abs{w}^2+\lambda^5\phi^5\abs{w_x}^2+\lambda^3\phi^3\abs{w_{xx}}^2+\lambda\phi\abs{w_{xxx}}^2\right)  \,dxdt  & \nonumber \\
+  \| P_1w \|_{L^2(Q)}^2   +  \| P_2w  \|_{L^2(Q)}^2 
&  \leq & C  \| Pw  \|_{L^2(Q)}^2  -   I_x .
\end{eqnarray*}
\end{lemm}

\proof

From hypothesis \eqref{hip1B} and the inequalities listed in \eqref{desigphi1}, we know that there exists $\delta >0$ such that
\begin{equation} \label{Normw}
  \sum_{k=0}^3 I(w_{kx}) \geq  \delta  \left\| w \right\|_{_{\la, \phi}}^2
\end{equation}
for a parameter $\lambda$ large enough.

Besides, from the definition \eqref{defR}, the fact   $\|q_i\|_{L^\infty((0,T) \times (0,1))}\leq m$ for $i=0,1,2$, and  \eqref{Phi1},   it is trivial to check that
\begin{equation}  \label{laR}
\begin{array}{rl}
\dis{ \|Rw\|^2_{L^2((0,T)\times(0,1))} }
& \dis{ \le C  \left( \la^6 \iint_Q \phi^6 |w|^2 + \la^2 \iint_Q \phi^2 |w_x|^2 + \la^2 \iint_Q \phi^2 |w_{xx}|^2         \right) } \vspace{0.3cm} \\
 & \dis{ \le C \la^{-1}   \left\| w \right\|_{_{\la, \phi}}^2}.    
\end{array}
\end{equation}

Thus, for $\lambda$ large enough, we have 
\begin{equation}  \label{P-R}
\begin{array}{rcl}
\dis{ \|P_1w\|_{L^2}^2 +  2\left< P_1w, P_2w\right> + \|P_2w\|_{L^2}^2 \,}
&  = & \dis{ \|Pw  - Rw\|^2_{L^2} } \vspace{0.3cm}  \\
& \dis{ \leq } & \dis{ 2 \left \|Pw \right\|^2_{L^2} + 2 \left \| Rw \right\|^2_{L^2} } \vspace{0.3cm}  \\
& \dis{ \leq } & \dis{ 2 \left \|Pw \right\|^2_{L^2} + C \la^{-1}   \left\| w \right\|_{_{\la, \phi}}^2  }.
\end{array}
\end{equation}

From Lemma \ref{Lemma_prodL2} and estimates (\ref{P-R}) and \eqref{Normw}, we conclude the proof of Lemma \ref{Carlemanw}.  \endproof

To complete the proof of Theorem \ref{CarlemanKS}, we have to deal with the norms fo $P_1 w$ and $P_2w$ appearing in Lemma \ref{Carlemanw}. From the definition of $P_2w$, 
and because (\ref{hip1B}) holds, we have 
\begin{equation*}  \label{P1P2-R}
 \frac{1}{\la \phi}|w_t|^2 
  \leq  \frac{2}{\la \phi} \left|P_2w \right|^2 + C  \left( \la^5  \phi^5 |w|^2 + \la^5  \phi^5|w_x|^2 + \la  \phi |w_{xxx}|^2         \right)   
\end{equation*}
and 
\begin{equation*}  \label{P1P2-R2}
 \iint_Q \frac{1}{\la \phi}|w_t|^2 
  \leq    C \iint_Q  \left|P_2w \right|^2 + C  \left\| w  \right\|_{\la,\phi}^2         
  \end{equation*}  
for $\lambda$ large enough. A similar result is proven for $(\sigma w_{xx})_{xx}$ and $P_1w$, and we then have
\begin{equation}  \label{P1P2-R12}
 \iint_Q \frac{1}{\la \phi} \left( |w_t|^2 + |(\sigma w_{xx})_{xx}|^2 \right)
  \leq    C \iint_Q  \left( \left|P_1w \right|^2 +  \left|P_2w \right|^2 \right) + C  \left\| w  \right\|_{\la,\phi}^2.        
  \end{equation}  
From (\ref{P1P2-R12}) and  Lemma \ref{Carlemanw} 
we obtain
\begin{multline}\label{Carl}
 \iint_Q \frac{1}{\lambda\phi}(\abs{w_t}^2+\abs{(\sig w_{xx})_{xx}}^2)+\lambda^7\phi^7\abs{w}^2+\lambda^5\phi^5\abs{w_x}^2+\lambda^3\phi^3\abs{w_{xx}}^2+\lambda\phi\abs{w_{xxx}}^2\\
\leq C \iint_Q \abs{Pw}^2
 - CI_x.
\end{multline}
To handle the terms in $I_x$, we note that for any $x \in (0,1)$ and $\lambda$ large enough,
\begin{equation*}\label{Carl2}
- C  \lambda \int_0^T \phi_x(x,t)\sig(x)\sig_{xx}(x)\abs{w_{xx}(x,t)}^2 dt
      \leq  C \lambda^3 \int_0^T \phi_x(x,t)^3\sig(x)^2\abs{w_{xx}(x,t)}^2 dt.
\end{equation*}
Then      
\begin{equation}\label{Carl3}
- C I_x 
      \leq  C \lambda^3 \int_0^T \phi_x(0,t)^3\sig(0)^2\abs{w_{xx}(0,t)}^2 dt
      + C  \lambda \int_0^T  \phi_x(0,t)\sig(0)^2\abs{w_{xxx}(0,t)}^2 dt
\end{equation}
and from (\ref{Carl}) and  (\ref{Carl3}) we obtain
\begin{multline}\label{Carl4}
 \iint_Q \frac{1}{\lambda\phi}\left(\abs{w_t}^2+\abs{(\sig w_{xx})_{xx}}^2\right)+
  \|w\|_{\la,\phi} ^2 
\leq C \iint_Q \abs{Pw}^2 \\
 + C \lambda^3 \int_0^T \phi_x(0,t)^3\sig(0)^2\abs{w_{xx}(0,t)}^2 dt
      + C  \lambda \int_0^T  \phi_x(0,t)\sig(0)^2\abs{w_{xxx}(0,t)}^2 dt. 
\end{multline}
Computing the derivatives of $e^{\la \phi}w$ it is trivial to prove that 
\begin{equation*}
\left| \partial_x^k  v \right|^2  = 
 \left| \partial_x^k ( e^{\la \phi } w ) \right|^2 \leq  C \sum_{j=0}^k \left|  \la^{k-j} \phi^{k-j} \partial_x^jw \right|^2
\end{equation*}
for each $k=0, \ldots, 3$.
Therefore
 \begin{multline*}
\int_0^T\int_0^1 e^{-2\la \phi} \left( \la^7 \phi^7 |e^{\la \phi}w|^2  + \la^5 \phi^5 |(e^{\la \phi}w)_x|^2   + \la^3 \phi^3   |(e^{\la \phi}w)_{xx}|^2
 + \la  \phi |(e^{\la \phi}w)_{xxx}|^2 \right) \,dxdt\\
 \leq  
 C \|w\|_{\la,\phi}.
 \end{multline*}
Considering finally that  $P w  = e^{-\la \varphi} Lv$, we obtain Carleman estimate \eqref{Carleman}.
\endproof

\begin{rem} 
We considered the function $\beta$ to be increasing. This allows the Carleman inequality 
to be obtained with boundary terms at $x=0$. If a decreasing function $\beta$ was used instead, then an inequality with boundary terms at $x=1$ would have been obtained. As 
discussed in the following section, the boundary terms in the Carleman inequality are related to the location of the observations in the inverse problem.
\end{rem}

\section{Inverse Problem}\label{InversePb}

In this section, the local stability of the nonlinear inverse problem stated in Theorem~\ref{PIKS} will be proved following the ideas of \cite{81bukhgeim-klibanov} and \cite{91klibanov-malinsky}. The proof is splited in several steps.\\

\noindent \textbf{Step 1. Local study of the inverse problem}\\
Let $\gamma$, $ \tilde \gamma$,  $y$ and  $\tilde y$ be defined as in Theorem \ref{PIKS}. If  we set  $u = y - \tilde y$ and $f =\tilde\gamma - \gamma$, 
then $u$ solves the following K-S equation:
\begin{equation}\label{KS2}
\left\{\begin{array}{ll}
u_t+(\sigma(x) u_{xx})_{xx}+  \gamma u_{xx}+  \tilde y u_x +  \tilde y_x u +uu_x = f(x)\tilde y_{xx}(x,t), \quad &\forall (t,x)\in Q,\\
u(t,0)=u(t,1)=0,\quad &\forall t\in(0,T),\\
u_x(t,0)=u_x(t,1)=0,\quad &\forall t\in(0,T),\\
u(0,x)=0,\quad &\forall x\in(0,1).
\end{array}\right.
\end{equation}

Then, in order to prove the stability of the inverse problem mentioned in the introduction, it is sufficient to obtain an estimate of $f$ in terms of $u_{xx}(\cdot, 0)$, $u_{xxx}(\cdot,0)$ and  $u(T_0, \cdot)$, where  $\tilde\gamma$ and $\tilde y$ are given, $\gamma\in  H^4(0,1)$ and $u$ is the solution of Equation~(\ref{KS2}).

We begin by deriving Equation (\ref{KS2}) with respect to time. Thus,  $v = u_t$ satisfies the following equation:
\begin{equation}\label{KS3}
\left\{\begin{array}{ll}
v_t+(\sigma v_{xx})_{xx}+ \gamma v_{xx}+ \tilde y v_x +  \tilde y_x  v =f\tilde y_{xxt} - g, \quad &\forall (t,x)\in Q,\\
v(t,0)=v(t,1)=0,\quad &\forall t\in(0,T),\\
v_x(t,0)=v_x(t,1)=0,\quad &\forall t\in(0,T),\\
v(0,x)=fR(x,0), \quad &\forall x \in(0,1),
\end{array}\right.
\end{equation}
where $g(x,t)  = u(x,t) y_{xt}(x,t) + u_x(x,t) y_t(x,t)$.\\

The proof of Theorem~\ref{PIKS} relies on the use of the Carleman estimate given in Theorem~\ref{CarlemanKS}. This result will be used twice. First, Equation (\ref{KS3}) allows to  estimate  $v$ in terms of $f$, $\tilde y_{xx}$ and $g$. Then, Equation (\ref{KS2}) will be used to handle the terms $u$  and $u_x$, which appear in the expression of the source term $g$. The details are given in the next step below.\\

\noindent \textbf{Step 2. First use of the Carleman estimate}\\
Similarly to the  proof of the Carleman estimate, we set $w = e^{-\lambda\phi}v$. Then, we work on the term
$$I = 2\int_0^1\int_0^{T_0} w(t,x) w_t (t,x)\,dt dx.$$

On the one hand, we can calculate $I$ and bound it from below. Indeed, 
using $w(0,x) = e^{-\lambda\phi(0,x)}v(0,x) = 0$ for all $x\in(0,1)$ and Equation (\ref{KS2}), we can easily obtain 
\begin{eqnarray*}
I  &=&  \int_0^1 \left|w\left( T_0,x\right)\right|^2\, dx \\
 & = &   \int_0^1e^{-2\lambda\phi\left( T_0,x\right)} \left| \left(f\tilde y_{xx} - (\sigma u_{xx})_{xx} - \gamma u_{xx} - \tilde y u_x -\tilde y_x u - uu_x\right)  \left( T_0, x\right)\right|^2 \, dx \\
 &\geq &   \int_0^1e^{-2\lambda\phi\left( T_0,x\right)} \abs{f(x)}^2 \left|\tilde y_{xx}\left( T_0, x\right)\right|^2 \, dx 
 - C \left\|u\left( T_0\right)\right\|_{H^4(0,1)}^2 
 - C \left\|u\left( T_0\right)\right\|_{H^1(0,1)}^4 
\end{eqnarray*}
where $C$ depends on  $\|\gamma\|_{L^\infty(0,1)}$,  $\| \tilde y(T_0) \|_{W^{1,\infty}(0,1)}$ and $\| \sigma  \|_{W^{2,\infty}(0,1)}$.

On the other hand, in order to estimate $I$ from above we apply the Carleman estimate~(\ref{Carl4}) 
to Equation~(\ref{KS3}) using $q_0=\tilde y_x$ and  $q_1=\tilde y$, which are uniformly  bounded in $L^\infty((0,T) \times (0,1))$ by the hypothesis in Theorem~\ref{PIKS}. We obtain that
\begin{eqnarray*}
I  &=&  2 \int_0^1\int_0^{T_0} w(t,x) w_t (t,x)\,dt dx\\
&\leq& \left(\int_0^1\int_0^{T_0}\lambda \phi(t,x) |w(t,x)|^2 \,dt dx \right)^{\frac 12} 
 \left( \int_0^1\int_0^{T_0}\frac{|w_t (t,x)|^2}{\lambda \phi(t,x)}\,dt dx\right)^{\frac 12}\\
&\leq& C \lambda^{-3} \int_0^1\int_0^{T} e^{-2\lambda\phi}\abs{f(x)\tilde y_{xxt}(x,t)}^2\,dxdt + C \lambda^{-3} \int_0^1\int_0^{T} e^{-2\lambda\phi}\abs{g(x,t)}^2\,dxdt\\
&&+ C \lambda^{-3}\int_0^T e^{-2\lambda\phi(0,t)}(\lambda^3\phi_x^3(0,t)\sig^2(0)\abs{v_{xx}(0,t)}^2+\lambda\phi_x(0,t)\sig^2(0)\abs{v_{xxx}(0,t)}^2)\,dt.
 \end{eqnarray*}

\noindent \textbf{Step 3. Second use of the Carleman estimate}\\
Considering that  $g  =  u y_{xt} + u_x y_t$, we will now use a Carleman estimate for the solution of Equation~(\ref{KS2}) in order to manage the term in $g$ of the previous inequality. The unknown trajectory $ y$ is nevertheless such that $ y_{xt}$ and $ y_{t}$ belong to $L^\infty(0,T;L^\infty(0,1))$ since 
$y\in \mathcal Z$. Thus, we have
\begin{eqnarray*}
\iint_Q e^{-2\lambda\phi}\abs{g}^2 
&\leq &2 \iint_Q e^{-2\lambda\phi}\abs{u}^2\abs{ y_{xt}}^2
+  2\iint_Q e^{-2\lambda\phi}\abs{u_x}^2\abs{ y_t}^2\\
&\leq & C \iint_Q e^{-2\lambda\phi}\left(\abs{u}^2 + \abs{u_x}^2\right).
 \end{eqnarray*}

Then we can apply Carleman estimate~(\ref{Carl4}) to equation (\ref{KS2}), using the identity $\tilde y_xu+uu_x=uy_x$, and taking $q_0=y_x$ and  $q_1= \tilde y$, which are  bounded in $L^\infty((0,T) \times (0,1))$. We can choose $\lambda_0$ as large as possible in Theorem~\ref{CarlemanKS}: we then obtain 
\begin{eqnarray*}
&&\iint_Q e^{-2\lambda\phi}\abs{g}^2 ~\leq  ~ C\lambda^{-5} \iint_Q e^{-2\lambda\phi}\left( \lambda^{7} \abs{u}^2 +\lambda^{5} \abs{u_x}^2\right)\\
&\leq& C \lambda^{-5} \iint_Q e^{-2\lambda\phi}\abs{f\tilde y_{xx}}^2 \\
&&~+ C \lambda^{-5}\int_0^T e^{-2\lambda\phi(0,t)}(\lambda^3\phi_x^3(0,t)\sig^2(0)\abs{u_{xx}(0,t)}^2+\lambda\phi_x(0,t)\sig^2(0)\abs{u_{xxx}(0,t)}^2)\,dt.
 \end{eqnarray*}
 
Gathering all the estimates of $I$  and $g$ that were obtained above, we have 
\begin{eqnarray*}
&&  \int_0^1e^{-2\lambda\phi\left( T_0,x\right)} \abs{f(x)}^2 \left|\tilde y_{xx}\left( T_0, x\right)\right|^2 \, dx 
 - C \left\|u\left( T_0\right)\right\|^2_{H^4(0,1)}
  - C \left\|u\left( T_0\right)\right\|_{H^1(0,1)}^4 
 \\
&\leq& C \lambda^{-3} \iint_Q e^{-2\lambda\phi}\abs{f\tilde y_{xxt}}^2  
+ C \lambda^{-8} \iint_Q e^{-2\lambda\phi}\abs{f\tilde y_{xx}}^2 \\
&&+~ C \lambda^{-8}\int_0^T e^{-2\lambda\phi(0,t)}(\lambda^3\phi_x^3(0,t)\sig^2(0)\abs{u_{xx}(0,t)}^2+\lambda\phi_x(0,t)\sig^2(0)\abs{u_{xxx}(0,t)}^2)\,dt\\
&&+~ C \lambda^{-3}\int_0^T e^{-2\lambda\phi(0,t)}(\lambda^3\phi_x^3(0,t)\sig^2(0)\abs{v_{xx}(0,t)}^2+\lambda\phi_x(0,t)\sig^2(0)\abs{v_{xxx}(0,t)}^2)\,dt.
 \end{eqnarray*}
From  the hypothesis of the theorem, we have 
$\tilde y \in C([0,T];H^{6}(0,1))$, $\tilde y_{t}\in C([0,T];H^{2}(0,1))$, and  $|\tilde y_{xx}(T_0,\cdot)|>\eta >0$  in $(0,1)$. Using also that the Carleman weight function satisfies \eqref{hip1P} thus $e^{-2\lambda\phi(t,x)}\leq e^{-2\lambda\phi(T_0,x)}$ in $(0,T)\times(0,1)$, we obtain
 \begin{eqnarray*}
&&\int_0^1e^{-2\lambda\phi(T_0,x)} \abs{f(x)}^2 \, dx \\
&&\leq C\Big( \lambda^{-3}  \int_0^1 e^{-2\lambda\phi(T_0,x)} \abs{f(x)}^2\,dx
+ \left\|u\left( T_0\right)\right\|^2_{H^4(0,1)} 
+  \left\|u\left( T_0\right)\right\|_{H^1(0,1)}^4 \\
&&+  \lambda^{-8}\int_0^T e^{-2\lambda\phi(0,t)}(\lambda^3\phi_x^3(0,t)\sig^2(0)\abs{u_{xx}(0,t)}^2+\lambda\phi_x(0,t)\sig^2(0)\abs{u_{xxx}(0,t)}^2)\,dt\\
&&+  \lambda^{-3}\int_0^T e^{-2\lambda\phi(0,t)}(\lambda^3\phi_x^3(0,t)\sig^2(0)\abs{v_{xx}(0,t)}^2+\lambda\phi_x(0,t)\sig^2(0)\abs{v_{xxx}(0,t)}^2)\,dt\Big).
 \end{eqnarray*}
 
 Therefore, the regularity of $\phi$ (that come from the assumptions on $\beta$ and $\phi_0$) allows to prove that choosing $\lambda_0$ large enough, we obtain the existence of a constant $C$ that depends on $r,K, T,\lambda_0, m$ such that $\forall \lambda>\lambda_0$,
 \begin{multline*}
 \|f(x)\|_{L^2(0,1)}^2
\leq C\Big( \|u( T_0,\cdot)\|_{H^4(0,1)}^2 
+   \left\|u\left( T_0,\cdot\right)\right\|_{H^1(0,1)}^4 \\
+ \|u_{xx}(\cdot,0)\|_{H^1(0,T)}^2
+\|u_{xxx}(\cdot,0)\|_{H^1(0,T)}^2\Big).
\end{multline*}
This estimate leads to the local stability of the initial inverse problem since $f= \tilde \gamma - \gamma$ and $u = y -\tilde y$  and we have proved Theorem~\ref{PIKS}. \\

\noindent{\bf Acknowledgments:} This work began while L. Baudouin and E. Cr\'epeau were visiting the Universidad T\'ecnica Federico Santa Mar\'ia on the framework of the MathAmsud project CIP-PDE. 
This study was partially supported by Fondecyt \#11080130, Fondecyt \#11090161, ANR C-QUID and CISIFS, and CMM-Basal grants.

\bibliographystyle{amsplain}

\end{document}